\documentclass[11pt]{article}
\usepackage[dvips]{epsfig}
\usepackage{amsmath,amssymb,euscript,graphicx,amsfonts,verbatim}
\usepackage{enumerate,color}
\usepackage{graphicx}
\usepackage{hyperref}

\usepackage{graphicx}

\newtheorem{theorem}{Theorem}[section]
\newtheorem{proposition}[theorem]{Proposition}
\newtheorem{lemma}[theorem]{Lemma}
\newtheorem{corollary}[theorem]{Corollary}

\newtheorem{conjecture}[theorem]{Conjecture}

\newtheorem{problem}[theorem]{Problem}
\newtheorem{example}[theorem]{Example}

\newtheorem{remark}[theorem]{Remark}

\newenvironment{proof}{{\noindent \sc Proof.}}{\hfill $\Qed$\\}

\newcommand{\Qed}{\rule{2.5mm}{3mm}}

\newcommand{\Cay}{\hbox{{\rm Cay}}}

\newcommand{\F}{{\cal{F}}}

\newcounter{case}

\renewcommand{\thecase}{\arabic{case}}

\newcounter{subcase}

\numberwithin{subcase}{case}

\begin{document}


\begin{center}
{\bf\large On maximum intersecting sets in direct and wreath product of groups} \\ [+4ex]
Ademir Hujdurovi\'c{\small$^{a,b,}$}\footnotemark,  
Klavdija Kutnar{\small$^{a,b,}$}\footnotemark$^{,*}$,  
\addtocounter{footnote}{0} 
Dragan Maru\v si\v c{\small$^{a, b, c,}$}\footnotemark 
\ and
\v Stefko Miklavi\v c{\small$^{a, b, c,}$}\footnotemark 
\\ [+2ex]
{\it \small 
$^a$University of Primorska, UP IAM, Muzejski trg 2, 6000 Koper, Slovenia\\
$^b$University of Primorska, UP FAMNIT, Glagolja\v ska 8, 6000 Koper, Slovenia\\
$^c$IMFM, Jadranska 19, 1000 Ljubljana, Slovenia}
\end{center}

\addtocounter{footnote}{-3}
\footnotetext{The work of Ademir Hujdurovi\'c  is supported in part by the Slovenian Research Agency (research program P1-0404 and research projects
J1-9110, N1-0102, J1-1691, J1-1694, J1-1695, N1-0140, 
N1-0159, J1-2451 and N1-0208).}
\addtocounter{footnote}{1}
\footnotetext{The work of Klavdija Kutnar  is supported in part by the Slovenian Research Agency (research program P1-0285 and research projects
 J1-9110, J1-9186, J1-1695, J1-1715, N1-0140, J1-2451,  J1-2481, and N1-0209).}
\addtocounter{footnote}{1}
\footnotetext{The work of Dragan Maru\v si\v c is supported in part by the Slovenian Research Agency (I0-0035, research program P1-0285
and research projects  J1-9108, 
J1-1694, J1-1695, N1-0140 and J1-2451).}
\addtocounter{footnote}{1}
\footnotetext{
The work of \v Stefko Miklavi\v c is supported in part by the Slovenian Research Agency (research program P1-0285
and research projects   J1-9110, 
J1-1695, N1-0140, N1-0159, J1-2451, and N1-0208).

~*Corresponding author e-mail:~klavdija.kutnar@upr.si}

\begin{abstract}
For a  permutation group $G$  acting on a set $V$,
a subset $I$ of $G$ is said to be an {\em intersecting set}
if for every pair of elements $g,h\in I$ there exists $v \in V$ such that
$g(v) = h(v)$.
The {\em intersection density} $\rho(G)$ of a transitive permutation
group $G$ is the maximum value of the quotient $|I|/|G_v|$ 
where $G_v$ is a stabilizer of a point $v\in V$ and 
$I$ runs over all intersecting sets in $G$. 
 If $G_v$ is the largest intersecting set in $G$ then
$G$ is said to have  the {\em Erd\H{o}s-Ko-Rado (EKR)-property}, and  moreover, 
$G$ has the {\em strict-EKR-property} if every intersecting set  of maximum
size in $G$ is a coset of a point stabilizer. 
Intersecting sets in $G$ coincide with  independent sets in the 
so-called
 {\em derangement graph} $\Gamma_G$, defined as  
the Cayley graph on $G$ with connection set consisting
of all derangements, that is,  fixed-point free elements of $G$.
In this paper a conjecture  regarding the existence of transitive permutation 
groups 
whose  derangement graphs are complete multipartite graphs, posed by
Meagher, Razafimahatratra and Spiga  in [{\em J.~Combin. Theory Ser. A} {\bf 180} (2021), 105390], is proved. The proof uses  direct product of groups. 
Questions regarding maximum intersecting sets in 
direct and wreath products of groups and the 
(strict)-EKR-property of these group products are also investigated.
In addition, some errors appearing in the literature on this topic are corrected.
\end{abstract}

\begin{quotation}
\end{quotation}

\begin{quotation}
\noindent 
{\em Math. Subj. Class.:} 05C25, 20B25.
\end{quotation}


\section{Introduction}
\label{sec:intro}
\noindent

The Erd\H{o}s-Ko-Rado theorem \cite{EKR} is one of the central result in extremal combinatorics. It gives a bound on the size of a family of intersecting $k$-subsets of a set and classifies the families satisfying the bound. Besides many
interesting proofs of this theorem, it has also been extended in various directions. For example,
Hsieh \cite{Hsieh} and Frankl and Wilson \cite{FW86} give a version for intersecting subspaces of a vector space over a
finite field, Rands \cite{Rands} extends it to intersecting
blocks in a design and Meagher and Moura \cite{MeagherM} prove a version for partitions. 

In this paper, we focus on the extension of the Erd\H{o}s-Ko-Rado theorem to permutation groups. For a finite set $V$ let $Sym(V)$ denote the corresponding symmetric group, and if $|V|=n$, we will use the notation $S_n$.
Given a permutation group $G\leq Sym(V)$, a subset $I$ of $G$ is called \emph{intersecting} if, for any two permutations $g$ and $h$ in $I$, there exists $v\in V$ such that $g(v)=h(v)$. 

The problem of determining the intersecting sets of 
a permutation group G can be formulated in
a graph-theoretic terminology. We denote by $\Gamma_G$ the \emph{derangement graph} of $G$, a Cayley graph 
whose  vertex set coincides
with $G$   and whose edges are the unordered pairs $\{g, h\}$ 
such that $gh^{-1}$ is a {\em derangement},
that is, $gh^{-1}$ is fixed-point-free. 
Now, an intersecting set of $G$ is simply an independent set or a coclique of $\Gamma_G$,
and similarly the classical Erd\H{o}s-Ko-Rado theorem translates into a classification of independent
sets of maximal cardinality of the Kneser graphs.

We say the group G has the EKR property if the size of any intersecting subset
of G is bounded above by the size of the largest point-stabilizer in G.

A permutation group $G$ is said to have the \emph{Erd\H{o}s-Ko-Rado property} (in short {\em EKR-property}), if the size of a maximum intersecting set is equal to the order of the largest point stabilizer, and is said to have 
the \emph{strict Erd\H{o}s-Ko-Rado property} (in short {\em strict-EKR-property})
 if every maximum intersecting set of $G$ is a coset of a point stabilizer. 
It is clear that strict-EKR-property implies EKR-property, but the converse does not hold. 
For example, the alternating group $A_4$ in its natural action 
on $\{1,2,3,4\}$ has the EKR-property, but not the strict-EKR-property. Namely, 
 $\{id,(1\,3\,2),(1\,4\,2)\}$ is a maximum intersecting set not coinciding with
a coset of a point stabilizer in $A_4$. In fact, for the natural action of the alternating groups, $A_4$ is the only alternating group without the strict-EKR-property (see \cite{AhmadiM1,KuWong}).
Furthermore,  the symmetric group 
$S_n$, $n\geq 2$, in its natural action has the strict-EKR-property
(see \cite{Cameron,DezaFrankl,GM,LaroseM}), and 
$p$-groups, $p$ a prime, have the EKR-property (see \cite{HKMM21,LSP}).
As for  doubly transitive permutation groups it was proved in \cite{MeagherSTEJC}
that they have the EKR-property, but there are infinitely many examples of such groups without the strict-EKR-property. 

The other end of the spectrum is more interesting: there
exist transitive permutation groups having the size of a maximum intersecting set bigger than the order of a point stabilizer. For example, the alternating group $A_5$ acting on 2-subsets of $\{1,2,3,4,5\}$ admits intersecting sets 
twice  as big as is the order of a point stabilizer in this action.
In \cite{LSP} a measure of how far a permutation group is from having 
the EKR property is introduced via  the concept of the so-called
intersection density.
The \emph{intersection density} of a transitive permutation group $G$, denoted by $\rho(G)$, is defined by $\rho(G)=\frac{|I|}{|G_v|}$, where $I$ is an intersecting set in $G$ of a maximum size, and $G_v$ is a point stabilizer. 
Clearly, $\rho(G)\geq 1$ with the equality holding 
if and only if $G$ has the EKR-property.
As proved in \cite{LSP}
intersection densities of transitive permutation groups can be arbitrarily large,
that is,  for every $M>0$ there exists a transitive permutation group $G$, with $\rho(G) > M$. 
Intersection densities of 
 transitive permutation groups of degree $2p$, $p$ a prime, 
have been determined to be either $1$ or $2$
(see \cite{HKMM21,R21} for details). 

In \cite{MRS21} Meagher, Razafimahatratra, and Spiga proved that the derangement graph of an arbitrary transitive permutation group of degree at least 3 contains a triangle and among other 
posed the following conjecture which is the
main motivation for this paper.

\begin{conjecture}
\label{conj:complete multipartite} {\rm\cite[Conjecture~6.6(1)]{MRS21}}
If $n$ is even but not a power of 2, then there is a transitive permutation group $G$ of degree $n$ such that $\Gamma_G$ is a complete multipartite graph with $n/2$ parts.
\end{conjecture}

The paper is organized as follows. 
In Section~\ref{sec2} we gather some useful results needed 
later on in the paper.
In Section~\ref{sec3}  we prove Conjecture~\ref{conj:complete multipartite}
 (see Theorem~\ref{thm:proof of the conjecture}). 
As a byproduct 
a characterization of the derangement graphs of direct products of groups 
is obtained (see Theorem~\ref{thm:direct product derangement graph}).
In particular, we prove that for transitive permutation groups $G$ and $H$, the direct product  $G\times H$ (in the canonical action)  has the (strict)-EKR-property if and only if 
both $G$ and $H$ have the (strict)-EKR-property, and that $\rho(G\times H)=\rho(G)\rho(H)$.  This implies that
by taking direct products of sufficiently many transitive permutation groups 
that do not have the EKR property a construction
of transitive permutation groups with arbitrarily large intersection densities
can be obtained.

In Section~\ref{sec:wreath products}  derangement graphs of 
wreath products of transitive permutation groups are 
characterized (see Lemma~\ref{lem:Derangement graph wreath product general} and Proposition~\ref{cor:wreath derangement graph H regular}).
Next, given transitive permutation groups $G$ and $H$ 
with $G\wr H$ denoting the corresponding wreath product,
we show that
 $\rho(G)\leq \rho(G\wr H)\leq \rho(G)\rho(H)$   (see Proposition~\ref{pro:density wreath bound}). 
Also, we give sufficient conditions for $G\wr H$ to have the strict-EKR-property (see Proposition~\ref{prop:generalizaion of AhmadiM}) and necessary and sufficient conditions for $G\wr H$ with $H$ being regular to have the strict-EKR-property (see Proposition~\ref{prop:wreath H regular strict EKR}).  

Along the way
certain inconsistencies in the material on this topic (see \cite{AhmadiM}
and Remarks~\ref{remark1} 
and~\ref{rem:Ahmadi Meagher mistake S_m wr S_n})
are taken care of.


\section{Preliminaries}
\label{sec2}
\noindent

%
Let ${\cal D}$ be the set of all derangements 
of a permutation group $G$. Following \cite{MRS21} 
we define the {\em derangement graph} of $G$
to be the graph $\Gamma_G=\Cay(G,{\cal D})$ with vertex set $G$
and edge set consisting of all edges $\{g,h\}$ such that
$gh^{-1}\in {\cal D}$. Therefore $\Gamma_G$ is the Cayley graph of $G$
with connection set ${\cal D}$, which is a loop-less simple graph
since ${\cal D}$
does not contain the identity element of $G$ and  
${\cal D}$ is inverse-closed. In the terminology of 
derangement graphs an intersecting set of $G$ is an 
independent set or a coclique in $\Gamma_G$. 

The largest size of a clique, the so-called  \emph{clique number},
in a graph $X$ is denoted by $\omega(X)$, 
and   the largest size of an independent set, that is, 
the \emph{independence number},  in $X$
is denoted by $\alpha(X)$.
The following classical upper bound on the size of the largest coclique in vertex-transitive graphs turns out to be quite useful when considering  intersection densities of  permutation groups. Namely, 
the derangement graph $\Gamma_G$ of a permutation group $G$ 
is always vertex-transitive.

\begin{lemma}
\label{lem:alpha-omega bound}
{\rm \cite[Corollary 4]{Cameron}}
Let $\Gamma$ be a vertex-transitive graph. 
Then the largest coclique in $\Gamma$ is of size  $\alpha(\Gamma)$
bounded by
$$
\alpha(\Gamma)\le\frac{|V(\Gamma)|}{\omega(\Gamma)},
$$
where $\omega(\Gamma)$ is the size of a maximum clique in $\Gamma$.
\end{lemma}

The next result is usually referred to as the  
``No-Homomorphism Lemma".

\begin{lemma}
{\rm \cite[Theorem 2]{Albertson}}
\label{lem:no-homomorphism}
Let $X$ be a
graph and  $Y$ be a vertex-transitive graph admitting 
a homomorphism from $X$ to $Y$. Then
$$
\frac{|V(X)|}{\alpha(X)} \leq \frac{|V(Y)|}{\alpha(Y)}.
$$
\end{lemma}

For graphs $X$ and $Y$ the {\em strong product} $X\boxtimes Y$ is the graph with vertex set $V(X)\times V(Y)$ such that 
$(x_1,y_1)$ is adjacent to $(x_2,y_2)$ if and only if 
one of the following three conditions holds:
\begin{enumerate}[(i)]
\itemsep=0pt
\item $x_1=x_2$ and $\{y_1,y_2\}\in E(Y)$;
\item $y_1=y_2$ and $\{x_1,x_2\}\in E(X)$;
\item $\{x_1,x_2\}\in E(X)$ and $\{y_1,y_2\}\in E(Y)$.
\end{enumerate}

The following result regarding maximum cliques in the strong product 
of graphs will be used for determining the intersection densities of direct product of groups. We use the following notation: for graphs $X$, $Y$ and for a set $C$ of vertices of  $X\boxtimes Y$ we denote by $p_X(C)=\{ x \in V(X) \mid \exists y \in V(Y) \textrm{ such that } (x,y) \in C\}$ and $p_Y(C)=\{ y \in V(Y) \mid \exists x \in V(X) \textrm{ such that } (x,y) \in C\}$ the projection of $C$ onto $V(X)$ and $V(Y)$, respectively.

\begin{lemma}
{\rm\cite[Lemma 7.3]{HIK}}
\label{lem:cliques in strong product}
Let $X$ and $Y$ be graphs and $C$ a maximum clique in $X\boxtimes Y$. Then $C = p_X(C) \times p_Y(C)$,
where $p_X(C)$ and $p_Y(C)$ is,  respectively, a  maximum clique 
in $X$ and $Y$.
\end{lemma}

For graphs $X$ and $Y$, the {\em the direct product} (also called the {\em tensor product}) $X \times Y$ is 
the graph with vertex set $V(X) \times V(Y)$ such that $(x_1,y_1)$ and $(x_2,y_2)$ are adjacent if and only if $\{x_1,x_2\} \in E(X)$ and $\{y_1,y_2\}\in E(Y)$.

The following result on the independence number of the direct product of vertex-transitive graphs
together with  its direct consequence in Lemma~\ref{lem:indepence number direct product}
will be needed in  Section~\ref{sec3}.

\begin{lemma}
{\rm \cite[Corollary~1.5]{Zhang12}}
\label{lem:Zhang}
Let $X$ and $Y$ be vertex-transitive graphs. Then 
$$\alpha(X\times Y)=\max\{\alpha(X)|V(Y)|,\alpha(Y)|V(X)|\}.$$
\end{lemma}

\begin{lemma}\label{lem:indepence number direct product}
Let $X$ be a vertex-transitive graph. Then $\alpha(X\times \cdots \times X)=\alpha(X)|V(X)|^{n-1}$.
\end{lemma}

For graphs $X$ and $Y$, the {\em wreath product $X[Y]$  of $X$ by $Y$} (also called the {\em lexicographic product})  is the graph with vertex set $V(X) \times V(Y)$ such that $(x_1,y_1)$ and $(x_2,y_2)$ are adjacent if and only if $\{x_1,x_2\} \in E(X)$ or $x_1=x_2$ and $\{y_1,y_2\}\in E(Y)$. The following result will be useful later in the paper. 

\begin{lemma}
\label{lem:lp}
{\rm  \cite[Theorem 1]{GS}}
Let $X$ and $Y$ be graphs. Then 
$$\alpha(X[Y])=\alpha(X)\alpha(Y).$$
\end{lemma}

Finally, the next 
result will be used in our analysis of intersection densities
of group products.

\begin{lemma}
\label{lem:density is decreasing}
{\rm \cite[Lemma~6.5]{MRS21}}
If $H\le G$ are transitive permutation groups then
$\rho(G)\le \rho(H)$.
\end{lemma}

\section{Derangement graphs of direct product of groups}
\label{sec3}

\noindent

Given finite sets $V$ and $W$  let
$G\leq Sym(V)$ and $H\leq Sym(W)$. We will consider two actions of $G\times H$. Namely, $G\times H$ acts naturally on $V\times W$ by the rule 
$$
(g,h):(v,w)\mapsto (g(v),h(w)) \textrm{ for $g\in G$ and $h\in H$}.
$$ 
Following \cite{AhmadiM} this action will be called the \emph{external direct product}. 
The other action of $G\times H$  is the action on $V \cup W$
(where $V$ and $W$ are considered as disjoint sets regardless of their structure), defined 
by the rule  
$$
(g,h):x \mapsto g(x) \textrm{ if $x\in V$ and $(g,h):x \mapsto h(x)$
 if $x\in W$}.
$$ 
Again, following \cite{AhmadiM} this action will be called 
the \emph{internal direct product}. Observe that the external direct product is transitive if and only if two factor groups are transitive, while
the internal direct product is never transitive.

\subsection{External direct products}
\noindent

The derangement graphs of the external direct products of groups were studied in \cite{AhmadiM}. In \cite[Lemma 4.1]{AhmadiM} it is claimed that $\Gamma_{G\times H}\cong \overline{\overline{\Gamma_G}\times \overline{\Gamma_H}}$. This is correct only if the complement $\overline{X}$ of a graph $X$ is defined in such a way that a loop is added to each vertex without a loop in  $X$. 
Since this is not a standard definition of a graph complement, which is not explicitly stated 
in  \cite{AhmadiM}, we give below a different characterization of derangement graphs of external direct products that uses the standard definition of a graph complement.

\begin{theorem}\label{thm:direct product derangement graph}
For finite sets $V$ and $W$  let
$G\leq Sym(V)$ and $H\leq Sym(W)$. 
Then the derangement graph $\Gamma_{G\times H}$ of  the
external direct product $G\times H$ acting on $V\times W$ satisfies
$$
\overline{\Gamma_{G\times H}}\cong \overline{\Gamma_G} \boxtimes \overline{\Gamma_H}.
$$
\end{theorem}
\begin{proof}
Let $(g_1,h_1),(g_2,h_2)$ be two different elements in $G\times H$.
Then 
\begin{align*}
&\{(g_1,h_1),(g_2,h_2)\}\in E(\overline{\Gamma_{G\times H}}) \Leftrightarrow\\
&(\exists (v,w) \in V\times W) \textrm{ such that } (v,w) \textrm{ is fixed by } (g_1g_2^{-1},h_1h_2^{-1}) \Leftrightarrow\\
&(\exists (v,w) \in V\times W) \textrm{ such that } (g_1g_2^{-1}\in G_v \wedge h_1h_2^{-1} \in H_w) \Leftrightarrow\\
&(g_1=g_2 \vee g_1 \sim_{\overline{\Gamma_G}} g_2) \wedge (h_1=h_2 \vee h_1 \sim_{\overline{\Gamma_H}} h_2) \Leftrightarrow\\
&(g_1=g_2 \wedge h_1 \sim_{\overline{\Gamma_H}} h_2) \vee (h_1=h_2 \wedge g_1 \sim_{\overline{\Gamma_G}} g_2 ) \vee (g_1 \sim_{\overline{\Gamma_G}} g_2 \wedge h_1 \sim_{\overline{\Gamma_H}} h_2) \Leftrightarrow\\
&\{(g_1,h_1),(g_2,h_2)\}\in E(\overline{\Gamma_G} \boxtimes \overline{\Gamma_H}).
\end{align*}
\end{proof}

Observe that if $H$ is a regular group, then $\Gamma_H$ is a complete graph, 
and so $\overline{\Gamma_H}$ is the empty graph,
giving us two immediate corollaries of 
Theorem~\ref{thm:direct product derangement graph}.

\begin{corollary}\label{cor:complement H regular}
For finite sets $V$ and $W$  let
$G\leq Sym(V)$ and $H\leq Sym(W)$. 
If $H$ is regular, then $\overline{\Gamma_{G\times H}}$ is isomorphic to
a disjoint union of $|H|$ copies of $ \overline{\Gamma_G}$.
\end{corollary}

\begin{corollary}\label{cor:complete multipartite}
For finite sets $V$ and $W$  let
$G\leq Sym(V)$ and $H\leq Sym(W)$.
If $H$ is regular and $\Gamma_ G$ is a complete multipartite graph with $k$ parts, then $\Gamma_{G\times H}$ is a complete multipartite graph with $k|H|$ parts.
\end{corollary}
\begin{proof}
If $\Gamma_G$ is a complete multipartite graph with $k$ parts, then $\overline{\Gamma_G}$ is a disjoint union of $k$ cliques. 
By Corollary~\ref{cor:complement H regular} it follows that $\overline{\Gamma_{G\times H}}$ is isomorphic a disjoint union of $|H|$ copies of $ \overline{\Gamma_G}$. Hence $\overline{\Gamma_{G\times H}}$ is a disjoint union of $k|H|$ cliques, and therefore $\Gamma_{G\times H}$ is a complete multipartite graph with $k|H|$ parts.
\end{proof}

We are now ready to prove Conjecture~\ref{conj:complete multipartite}.

\begin{theorem}\label{thm:proof of the conjecture}
If $n$ is even but not a power of 2, then there is a transitive permutation group $G$ of degree $n$ such that $\Gamma_G$ is a complete multipartite graph with $n/2$ parts.
\end{theorem}
\begin{proof}
Let $n=2^ak$ with $k \ge 3$ odd and $a\geq 1$.
By \cite[Lemma 5.3]{MRS21} there exists a transitive permutation group $K$ of degree $2k$ such that $\Gamma_K$ is a complete multipartite graph with $k$ parts. Let $H$ be a regular group of degree $2^{a-1}$. Then by Corollary~\ref{cor:complete multipartite} it follows that $\Gamma_{K\times H}$ is a complete multipartite graph with $2^{a-1}k=n/2$ parts.
\end{proof}

Recall that maximum cocliques in $\Gamma_G$ (equivalently, maximum cliques of $\overline{\Gamma_G}$) correspond to maximum intersecting sets of $G$. 
Consequently, if $G$ is transitive, $\rho(G)=\omega(\overline{\Gamma_G})/|G_v|$.

\begin{theorem}\label{thm:direct-intersection density}
Let $G\leq Sym(V)$ and $H\leq Sym(W)$ be transitive permutation groups
and let $G\times H$ be the external direct product acting on $V\times W$. 
Then 
\begin{enumerate}[(i)]
\itemsep=0pt
\item $\rho(G\times H)=\rho(G)\rho(H)$, and
\item   $G\times H$ has 
the strict-EKR-property if and only if both $G$ and $H$ have the 
strict-EKR-property.
\end{enumerate}
\end{theorem}

\begin{proof}
Recall that $\rho(G\times H)=\omega(\overline{\Gamma_{G \times H}})/|(G\times H)_{(v,w)}|$ where   $(v,w)\in V\times W$. 
Since $\overline{\Gamma_{G\times H}}\cong \overline{\Gamma_G} \boxtimes \overline{\Gamma_H}$, by Lemma~\ref{lem:cliques in strong product} it follows that $\omega(\overline{\Gamma_{G \times H}})=\omega (\overline{\Gamma_{G}}) \omega (\overline{\Gamma_{H}})$. 
Observe that $(G\times H)_{(v,w)}=G_v\times H_w$. It follows that
$\rho(G\times H)=\omega(\overline{\Gamma_{G \times H}})/|(G\times H)_{(v,w)}|=\omega (\overline{\Gamma_{G}}) \omega (\overline{\Gamma_{H}})/|G_v||H_v|=\rho(G)\rho(H)$, proving (i).

As for (ii),
let $\F$ be a maximum intersecting set of $G\times H$ containing $(1_G,1_H)$ and let $Q$ be the corresponding maximum independent set in $\Gamma_{G\times H}$, or equivalently a maximum clique in $\overline{\Gamma_{G\times H}}$.
By Theorem~\ref{thm:direct product derangement graph} 
we have $\overline{\Gamma_{G\times H}}\cong \overline{\Gamma_G} \boxtimes \overline{\Gamma_H}$. 
By Lemma~\ref{lem:cliques in strong product} it follows that $Q=p_{\overline{\Gamma_G}}(Q)\times p_{\overline{\Gamma_H}}(Q)$
where $p_{\overline{\Gamma_G}}(Q)$ and $p_{\overline{\Gamma_H}}(Q)$
are maximum cliques of $\overline{\Gamma_G}$ and $\overline{\Gamma_H}$ containing $1_G$ and $1_H$, respectively. 
Now, if $G$ and $H$ have the strict-EKR-property, then
 $p_{\overline{\Gamma_G}}(Q)$ and $p_{\overline{\Gamma_H}}(Q)$
(viewed as subsets of $G$ and $H$ respectively) 
are vertex stabilizers in $G$ and $H$, and consequently $Q$ 
(and equivalently $\F$) is a vertex-stabilizer in $G\times H$.
Conversely, if $Q$ is a vertex-stabilizer in $G\times H$, it follows that $p_{\overline{\Gamma_G}}(Q)$ and $p_{\overline{\Gamma_H}}(Q)$ are vertex-stabilizers in $G$ and $H$, respectively.
\end{proof}

\begin{remark}\label{remark1}
{\rm
In \cite[Theorem 4.3]{AhmadiM} it is claimed that the external direct product of groups has the (strict)-EKR-property if and only if each of the factors has the (strict)-EKR-property, a claim made in Theorem~\ref{thm:direct-intersection density} too. We would like to note, however, that the proof of \cite[Theorem 4.3]{AhmadiM} is incomplete, as it relies on \cite[Lemma 4.1]{AhmadiM} (which is incorrect unless the complement of a graph is defined to include loops) and \cite[Lemma 4.2]{AhmadiM}, the proof of which seems to be missing some details.
}
\end{remark}

\begin{remark}{\rm
Given a group $G$ with $\rho(G)>1$, one can construct a group 
with arbitrarily large intersection density
by taking sufficiently many copies of $G$ in
the direct product $G\times G\times \cdots \times G$.}
\end{remark}

\subsection{Internal direct products}
\noindent

The study of intersecting sets in internal direct products was initiated by Ku and Wong in \cite{KuWong}. They gave necessary and sufficient conditions for the direct product of symmetric groups to have the strict-EKR-property. In \cite{AhmadiM} the derangement graphs of internal direct products were characterized.

\begin{lemma}
{\rm \cite{AhmadiM}}
\label{lem:internal direct derangement}
Let $V_1,\ldots,V_n$ be disjoint sets and  $G_i\leq Sym(V_i)$ $i\in \{1,\ldots,n\}$. Let $G=G_1\times \cdots \times G_n$ be the internal direct product acting on $V_1\cup \cdots \cup V_n$. Then $\Gamma_G\cong \Gamma_{G_1} \times \cdots \times \Gamma_{G_n}$. Moreover, if all groups $G_i$ have the EKR-property, then $G$ also has the EKR-property.
\end{lemma}

No analogue to Lemma~\ref{lem:internal direct derangement} holds 
for the strict-EKR-property. Namely the internal direct product of groups
each having the strict-EKR-property does not need to have the strict EKR-property. For example, $S_3 \times S_3$ does not have the strict EKR-property.

We are  interested in determining conditions under which does the internal direct product of the form $G\times G \times \cdots \times G$ have the strict-EKR-property.
As we will see in Section~\ref{sec:wreath products} 
this is important for understanding when does the wreath product of groups have the strict-EKR-property.
With the exception of the internal direct products of full symmetric groups, not much has been done on this problem. 
This is not surprising, as on the one hand, these groups are not transitive (and restriction to investigate  transitive permutation groups is common), and on the other, the problem is quite hard. 
Namely, one would need to characterize maximum independent sets in the direct product of derangement graphs of the factor groups, which is not
an easy task in view of the fact that
the structure of maximum independent sets in direct products is 
usually not all that simple.  The situation changes somewhat with 
the so-called MIS-normal direct products introduced below.

The direct product $X_1\times \cdots \times X_n$ of graphs $X_1,\ldots , X_n$ is said to be \emph{MIS-normal} (maximum-independent-set-normal) if every maximum independent set of it is the preimage of an independent set of one factor under projections.
Following Zhang \cite{Zhang12} we say that a graph $\Gamma$ is \emph{IS-primitive} if there is no non-maximum independent set $A$ such that $|A|/|N[A]|=\alpha(\Gamma)/|V(\Gamma)|$,
where $N[A]$ denotes the closed neighbourhood of $A$. Note that a disconnected vertex-transitive graph is not IS-primitive.
The following result follows from \cite[Theorem 1.4]{Zhang12} and \cite[Corollary 3.2]{Zhang11}.
\begin{lemma}\label{lem:MIS normal}
Let $X$ be a non-bipartite vertex-transitive graph. Then   $X\times \cdots \times X$ is MIS-normal if and only if $X$ is IS-primitive.
\end{lemma}

\begin{example} {\rm
Let $G$ be a regular group. Then $\Gamma_G$ is IS-primitive. Namely, the graph $\Gamma_G$ is a complete graph, and hence maximum independent sets are of size 1.}
\end{example}

\begin{lemma}\label{lem:internal direct strict EKR main}
Let $G\leq Sym(V)$ be a transitive permutation group with $|V|\geq 3$. 
Then the internal direct product $G\times \cdots \times G$ has the strict-EKR-property if and only if $G$ has the strict-EKR-property and $\Gamma_G$ is IS-primitive.
\end{lemma}

\begin{proof}
Let $G^n=G\times \cdots \times G$.
Suppose first that $G^n$ has the strict-EKR-property. If $G$ does not have the 
strict-EKR-property, then there exists a maximum intersecting set $S$ in  $G$ 
containing the identity which is not a point stabilizer. Then $I=S\times G 
\times \cdots \times G$ is a maximum intersecting set in $G^n$ which is not a 
point stabilizer, contradicting the assumption that $G^n$ has 
the strict-EKR-property.  Since $G$ is of degree at least 3, it follows by \cite[Theorem 1.5]{MRS21} that $\Gamma_G$ is non-bipartite.
Therefore, if $\Gamma_G$ is not IS-primitive, then 
Lemma~\ref{lem:MIS normal} implies that $\Gamma_{G^n}$ is not MIS-normal, that is, there exists 
a maximum independent set in $\Gamma_{G^n}$ which is not a preimage of a 
maximum independent set in $\Gamma_G$. As vertex stabilizers in $G^n$ have 
the form $G\times \cdots \times G \times G_v \times G \times \cdots \times G$, 
we conclude that $G^n$ does not have the strict-EKR-property, a contradiction
showing that $\Gamma_G$ is IS-primitive.

Conversely, suppose that $G$ has the strict-EKR-property 
and that $\Gamma_G$ is IS-primitive. 
Lemma~\ref{lem:internal direct derangement} implies
 that $\Gamma_{G^n}\cong \Gamma_G\times \cdots \times \Gamma_G$.
Since $G$ is of degree at least 3, it follows by \cite[Theorem 1.5]{MRS21} 
that $\Gamma_G$ is non-bipartite.
Since $\Gamma_G$ is IS-primitive, it follows by Lemma~\ref{lem:MIS normal} 
that $\Gamma_{G^n}$ is MIS-normal. Let $I$ be a maximum intersecting set in 
$G^n$ containing the identity. Then $I$ is a maximum independent set in $
\Gamma_{G^n}$, and since the graph is MIS-normal, we may conclude, 
 without loss 
of generality, that $I=I_0\times  G \times \cdots \times G$, where $I_0$ is a 
maximum independent set in $\Gamma_G$. Since $G$ has the strict-EKR-
property, it follows that $I_0$ is a stabilizer of a point in $G$. Hence $I$ is a 
stabilizer of a point in $G^n$, showing that $G^n$ has the strict-EKR-property.	
\end{proof}

We close the section with an open problem.

\begin{problem}
Characterize transitive permutation groups $G$ such that $\Gamma_G$ is IS-primitive.
\end{problem}

\section{Derangement graphs of wreath products}
\label{sec:wreath products}
\noindent

Let $G$ be a permutation group acting on a set $V$ and $H\le S_n$ a permutation group acting on the set $N=\{1,\ldots,n\}$. The {\it wreath product of $G$ by $H$} denoted by $G\wr H$ is the set of all permutations $((g_1,\ldots,g_n),h)$ of 
$V\times N$ (where $g_1,\ldots,g_n \in G$ and $h \in H$) such that $((g_1,\ldots,g_n),h): (a,i)\mapsto (g_i(a),h(i))$.


Observe that  for $((g_1,\ldots,g_n),h), ((g'_1,\ldots,g'_n),h') \in G\wr H$ we have 
$$((g_1,\ldots,g_n),h)\cdot((g'_1,\ldots,g'_n),h')=(g_{h'(1)}g'_1,\ldots,g_{h'(n)}g'_n),hh')$$
and
$$
((g_1,\ldots,g_n),h)^{-1}=((g_{h^{-1}(1)}^{-1},\ldots,g_{h^{-1}(n)}^{-1}),h^{-1}).
$$
Observe that the wreath product $G\wr H$ is isomorphic to 
the semidirect product $G^n\rtimes H$, where $G^n$ is the internal direct product of $G$ with $n$ factors.
The main goal of this section is to characterize derangement graphs of wreath products, and to study the intersection density, EKR-property and strict-EKR-property of wreath products. We start with the following lemma.

\begin{lemma}\label{lem:adjacency in derangement wreath}
\begin{sloppy}
Let $G\leq Sym(V)$ and $H\leq S_n$ be permutation groups. Vertices $((g_1,\ldots,g_n),h)$ and $((g'_1,\ldots,g'_n),h')$ are adjacent in $\Gamma_{G\wr H}$ if and only if $g_i$ is adjacent with $g'_i$ in $\Gamma_G$ for every $i \in \{1,\ldots,n\}$ such that $h(i)=h'(i)$.
\end{sloppy}
\end{lemma}
\begin{proof}

\begin{sloppy}
Clearly, two vertices $\mathbf{g}=((g_1,\ldots,g_n),h)$ and $\mathbf{g}'=((g'_1,\ldots,g'_n),h')$ are adjacent if and only if 
$\mathbf{g'}\mathbf{g}^{-1}=((g'_{h^{-1}(1)}g_{h^{-1}(1)}^{-1},\ldots,g'_{h^{-1}(n)}g_{h^{-1}(n)}^{-1}),h'h^{-1})$ is a derangement. 
\end{sloppy}
If $j\not \in fix(h'h^{-1})$ then $\mathbf{g'}\mathbf{g}^{-1}$ does not fix a
 point in $V\times \{j\}$ regardless of what are the values of $g_k$  and $g'_k$ for $k\in \{1,\ldots,n\}$. 

Suppose that $j \in fix(h'h^{-1})$. Observe that $\mathbf{g'}\mathbf{g}^{-1} (v,j)=(g'_{h^{-1}(j)}g_{h^{-1}(j)}^{-1}(v),j)$. Hence $\mathbf{g'}\mathbf{g}^{-1}$ does not fix a point in $V\times \{j\}$ if and only if $g'_{h^{-1}(j)}g_{h^{-1}(j)}^{-1}$ is a derangement. 
We conclude that $\mathbf{g'}\mathbf{g}^{-1}$ is a derangement if and only if $g'_{h^{-1}(j)}g_{h^{-1}(j)}^{-1}$ is a derangement for each $j \in fix(h'h^{-1})$. This is equivalent with $g'_ig_i^{-1}$ being a derangement for each $i\in \{1,\ldots,n\}$ such that $h(i) \in fix(h'h^{-1})$. 
The last condition is equivalent with $h(i)=(h'h^{-1})(h(i))=h'(i)$.
To summarize, $((g_1,\ldots,g_n),h)$ and $((g'_1,\ldots,g'_n),h')$ are adjacent in $\Gamma_{G\wr H}$ if and only if $g_i$ is adjacent with $g'_i$ in $\Gamma_G$ for every $i \in \{1,\ldots,n\}$ such that $h(i)=h'(i)$.
\end{proof}

Let $G\leq Sym(V)$ and $H\leq S_n$ be permutation groups, and let $h,h'\in H$. Denote with $\Gamma_{G\wr H}[h,h']$ the subgraph of $\Gamma_{G\wr H}$  with vertex set $G^n \times \{h,h'\}$ containing all the edges with one endvertex in $G^n \times \{h\}$ and the other in $G^n \times \{h'\}$. Observe that for $h=h'$, the graph $\Gamma_{G\wr H}[h,h]$ is the subgraph of $\Gamma_{G\wr H}$ induced by the set $G^n \times \{h\}$.
We will  denote it by $\Gamma_{G\wr H}[h]$. As for $h\neq h'$, the subgraph $\Gamma_{G\wr H}[h,h]$ is the bipartite subgraph of $\Gamma_{G\wr H}$ containing all the edges between sets $G^n \times \{h\}$ and $G^n \times \{h'\}$. 
We now  give a characterization of the derangement graph of a wreath product of groups. To simplify the notation, we denote by $K^*_n$ the complete graph of order $n$ with a loop at each vertex.

\begin{lemma}\label{lem:Derangement graph wreath product general}
Let $G\leq Sym(V)$ and $H\leq S_n$ be permutation groups and let $h,h'\in H$ be distinct. Then:
\begin{enumerate}[(i)]
\itemsep=0pt
\item $\Gamma_{G\wr H}[h] \cong \Gamma_G \times \Gamma_G \times \cdots \times \Gamma_G$; and
\item $\Gamma_{G\wr H}[h,h']\cong X_1 \times X_2 \times \cdots \times X_n \times K_2$, where  $X_i=\Gamma_G$ if $h(i) = h'(i)$ and $X_i=K^*_{|G|}$ if $h(i) \neq h'(i)$. 
\end{enumerate}
\end{lemma}
\begin{proof}
Applying Lemma~\ref{lem:adjacency in derangement wreath} for vertices $((g_1,\ldots,g_n),h)$ and $((g'_1,\ldots,g'_n),h)$ it follows that $((g_1,\ldots,g_n),h)$ and $((g'_1,\ldots,g'_n),h)$ are adjacent in $\Gamma_{G\wr H}$ if and only if $g_i$ is adjacent to $g'_i$ in $\Gamma_G$ for each $i\in \{1,\ldots,n\}$. 
Now the claim $(i)$ holds directly by the definition of direct product of graphs.  
Similarly, (ii) follows from 
Lemma~\ref{lem:adjacency in derangement wreath} and the definition of direct product of graphs.
\end{proof}

In the case when $H\leq S_n$ is regular, the derangement graph of $G\wr H$ has a particularly nice structure, given as follows.

\begin{proposition}\label{cor:wreath derangement graph H regular}
Let $G\leq Sym(V)$ and $H \leq S_n$. If $H$ is regular, then 
$
\Gamma_{G\wr H}\cong K_n [\Gamma_G \times \Gamma_G \times \cdots \times \Gamma_G].
$
\end{proposition}
\begin{proof}
Since $H$ is regular,  we have $h(i)\neq h'(i)$ 
for any two distinct elements $h,h'\in H$ and for every $i\in \{1,\ldots,n\}$.
Then Lemma~\ref{lem:Derangement graph wreath product general} 
implies that $\Gamma_{G\wr H}[h,h']\cong K^*_{|G|}\times \cdots \times K^*_{|G|}\times K_2$. Observe that the last graph is isomorphic to a complete bipartite graph of order $2|G|^n$, that is, there are all possible edges between sets $G^n\times \{h\}$ and $G^n\times \{h'\}$ in $\Gamma_{G\wr H}$. Since by Lemma~\ref{lem:Derangement graph wreath product general}(i) the subgraph induced by $G^n\times \{h\}$ is isomorphic to $\Gamma_G \times \Gamma_G \times \cdots \times \Gamma_G$, by the definition of wreath product of graphs, it follows that $\Gamma_{G\wr H}\cong K_n [\Gamma_G \times \Gamma_G \times \cdots \times \Gamma_G]$.
\end{proof}

The next proposition  gives bounds on the intersection density of the wreath product of groups.

\begin{proposition}\label{pro:density wreath bound}
Let $G\leq Sym(V)$ and $H\leq S_n$ be transitive permutation groups. Then  $\rho(G)\leq \rho(G\wr H)\leq \rho(G)\rho(H)$.
\end{proposition}
\begin{proof}
Let $I_G$ be a maximum intersecting set in $G$ and let $I=I_G \times G \times \cdots \times G \times H_1\subset G\wr H$. We claim that $I$ is an intersecting set in $G \wr H$. 
Let
$\mathbf{g}=((g_1,\ldots,g_n),h)$ and $\mathbf{g}'=((g'_1,\ldots,g'_n),h')$ be arbitrary elements of $I$. 
Since $I_G$ is an intersecting set in $G$ there exists $v\in V$ such that $g_1(v)=g'_1(v)$. 
It follows that $\mathbf{g}(v,1)=(g_1(v),1)=(g'_1(v),1)=\mathbf{g}'(v,1)$. We conclude that $I$ is indeed an intersecting set in $G\wr H$.

Observe that $|I|=|I_G||H_1||G|^{n-1}$ and $|(G\wr H)_{(v,i)}|=\frac{|G|^n |H|}{|V|\cdot n}$. Using the orbit stabilizer theorem for $G$ and $H$ it follows that $\frac{|I|}{|(G\wr H)_{(v,i)}|}=\frac{|I_G|}{|G_1|}=\rho(G)$. Since $I$ is an intersecting set in $G\wr H$ we have $\rho(G\wr H) \geq \frac{|I|}{|(G\wr H)_{(v,i)}|}=\rho(G)$.

Lemma~\ref{lem:Derangement graph wreath product general} implies that $\Gamma_H[\Gamma_G \times \cdots \times \Gamma_G]$ is a subgraph of $\Gamma_{G\wr H}$, and so there exists a homomorphism from $\Gamma_H[\Gamma_G \times \cdots \times \Gamma_G]$ to $\Gamma_{G\wr H}$. Using Lemmas~\ref{lem:no-homomorphism}, \ref{lem:indepence number direct product}  and \ref{lem:lp} it follows that
\begin{equation}\label{eq:no homomorphism}
\frac{|G|^n\cdot |H|}{\alpha(\Gamma_H)\cdot \alpha(\Gamma_G) \cdot |G|^{n-1}} \leq \frac{|G|^n\cdot |H|}{\alpha(\Gamma_{G\wr H})}.
\end{equation}
Equation~\eqref{eq:no homomorphism} implies that $\alpha(\Gamma_{G\wr H})\leq \alpha(\Gamma_G) \alpha(\Gamma_H) |G|^{n-1}$. It follows that 
$$
\rho(G\wr H) =\frac{\alpha(\Gamma_{G\wr H})}{|G_v||G|^{n-1}|H_1|} \leq
\frac{\alpha(\Gamma_H)\alpha(\Gamma_G)}{|G_v||H_1|}=\rho(G)\cdot \rho(H). 
$$

\end{proof}

\begin{corollary}
\label{cor:wreath1}
Let $G\leq Sym(V)$ be a transitive permutation group, and let $H\leq S_n$ be a transitive permutation group having the EKR-property. Then $\rho(G\wr H)=\rho(G)$.
\end{corollary}

Determining the exact value of $\rho(G\wr H)$ for all transitive permutation groups $G$ and $H$ is not an easy task. We conjecture the following is true.

\begin{conjecture}
Let $G\leq Sym(V)$ and $H\leq S_n$ be transitive permutation groups. Then $\rho(G\wr H)=\rho(G)$.
\end{conjecture}

Now we turn our attention to study the following problem.

\begin{problem}
	\label{prob1}
Let $G\leq Sym(V)$ and $H\leq S_n$.
Give necessary and sufficient conditions for $G\wr H$ to have the strict-EKR-property.
\end{problem}


If $H$ is regular, then $\Gamma_{G\wr H}$ is isomorphic to $K_n[\Gamma_G \times \cdots \times \Gamma_G]$, where $K_n$ denotes the complete graph on $n$ vertices.  Observe that every independent set of  $K_n[\Gamma_G \times \cdots \times \Gamma_G]$ is of the form $\{v\} \times I$, where $v$ is a vertex of $K_n$ and $I$ is an independent set of $\Gamma_G \times \cdots \times \Gamma_G$. This shows that in the case when $H$ is regular, Problem \ref{prob1} reduces to the question when does the internal direct product $G\times \cdots \times G$ have the strict-EKR-property, which 
was answered in Lemma~\ref{lem:internal direct strict EKR main}.

\begin{proposition}\label{prop:wreath H regular strict EKR}
Let $G\leq Sym(V)$ be transitive with $|V|\geq 3$ and let $H\leq S_n$ be regular. Then $G\wr H$ has the strict-EKR-property if and only if $G$ has the strict-EKR-property and $\Gamma_G$ is IS-primitive.
\end{proposition}

\begin{corollary}\label{cor:S_3 wr S_2}
The group $S_3\wr S_2$ does not have the strict-EKR-property.
\end{corollary}
\begin{proof}
Note that $S_2$ is regular. Observe that $\Gamma_{S_3}$ is disconnected (a disjoint union of two triangles), hence $\Gamma_{S_3}$ is not IS-primitive. The result now follows from Proposition~\ref{prop:wreath H regular strict EKR}.
\end{proof}

The next proposition generalizes \cite[Proposition 6.5]{AhmadiM}.

\begin{proposition}\label{prop:generalizaion of AhmadiM}
Let $G\leq Sym(V)$, $|V|\geq 3$, and $H\leq S_n$ be 
transitive permutation groups. 
If the internal direct product $G^n$ of $n$ factors of the group $G$
has the strict-EKR-property and $H$ has the EKR-property, 
then $G\wr H$ has the strict-EKR-property. 
\end{proposition}

\begin{proof}
Since $P=G^n$ has the strict-EKR-property, it follows by Lemma~\ref{lem:internal direct strict EKR main} that $G$ has the strict-EKR-property. Let $S$ be an intersecting set of $G\wr H$ of maximum size, that is $|S|=|G|^{n-1}|G_v||H_1|$. Without loss of generality we may assume that $S$ contains the identity of $G\wr H$. Since $\Gamma_{G\wr H}$ contains a subgraph isomorphic to $\Gamma_H[\Gamma_P]$ it follows that $S$ is also an independent set in $\Gamma_H[\Gamma_P]$. Observe that $\alpha(\Gamma_H[\Gamma_P])=\alpha(\Gamma_H)\alpha(\Gamma_P)=|S|$. Hence $S$ is an intersecting set of $\Gamma_H[\Gamma_P]$ of maximum size.

For $h\in H$ let $S_h=S \cap (G^n \times \{h\})$. Then for every $h\in H$, we have $S_h=\emptyset$ or $S_h$ is a maximum intersecting set in $P$. Since $P$ has the strict-EKR-property, it follows that 
$S_h\neq \emptyset$ implies that $S_h$ is a coset of a point stabilizer in $P$. Since $S$ contains the identity of $G\wr H$, it follows that $S_h=P_{(v(h),i(h))}$ for some $v(h) \in V$ and some $i(h)\in \{1,\ldots,n\}$. Let $h$ and $h'$ be two distinct elements of $H$ such that $S_h$ and $S_{h'}$ are non-empty. It follows that   $h'h^{-1}$ is not a derangement. 

Suppose that $i(h) \neq i(h')$. 
 Let $g\in G$ be a derangement. Consider elements $((g_1,\ldots,g_n),h)$ and 
$((g'_1,\ldots,g'_n),h')$ in $G\wr H$ where $g_j=g$ for $j\neq i(h)$ and 
$g_{i(h)}=1$, while $g'_j=1$ for $j\neq i(h)$ and $g'_{i(h)}=g$. Observe that $ 
g'_jg_j^{-1}$ is a derangement for each $j\in \{1,\ldots,n\}$ hence $((g_1,
\ldots,g_n),h)$ and $((g'_1,\ldots,g'_n),h')$ are adjacent in $\Gamma_{G\wr H}
$ which contradicts the fact that they are both contained in $S$.
Consequently,  
 $i(h)=i(h')$. 

To simplify the notation, let $i(h)=i(h')=k$.
Suppose  that $P_{(v(h),k)}\neq P_{(v(h'),k)}$, that is, $G_{v(h)}\neq G_{v(h')}$. 
Choose $f\in G_{v(h)}$ and $f'\in  G_{v(h')}$ in such a way that $f'f^{-1}$ is a 
derangement. Such elements always exist for otherwise $G_{v(h)} \cup 
G_{v(h')}$ would be an independent set in $\Gamma_G$ of size greater than 
the size of a point stabilizer. 
Consider elements $((g_1,\ldots,g_n),h)$ and $((g'_1,\ldots,g'_n),h')$ in $G\wr 
H$ where $g_j=g$  for $j\neq k$ and 
$g_{k}=f$, while $g'_j=1$ for $i\neq k$ and $g'_{k}=f'$. As in the 
previous paragraph it follows that $((g_1,\ldots,g_n),h)$ and $((g'_1,
\ldots,g'_n),h')$ are adjacent in $\Gamma_{G\wr H}$, contradicting the fact that they are both contained in $S$, an independent set in 
$\Gamma_{G\wr H}$. This shows that $P_{(v(h),k)}= 
P_{(v(h'),k})$. Let $w=v(h)$. It follows that $S=P_{(w,k)}\times S_H$ where $S_H
$ is a maximum intersecting set in $H$. Since $H$ has the EKR-property, it 
follows that $|S_H|=|H_k|$. Since $S$ contains the identity 
of $G\wr H$, it follows 
that $1_H\in S_H$.
 
Suppose that $S_H\neq H_k$. Then there exist $h,h'\in S_H$ 
such that $h(k)\neq h'(k)$. Namely, if that was not the case we would have 
that $h^{-1}h'(k)=k$ for every $h,h'\in S_H$. 
Then by taking $h'=1_H$ we would have $h(k)=k$ for every $h\in H_j$  implying that $S_H=H_k$, contrary to the assumption that $S_H\neq H_k$.
Let us now consider elements $((g_1,\ldots,g_n),h)$ 
and $((g'_1,\ldots,g'_n),h')$ in $G\wr H$ where 
$g_j=g$ 
for $j\neq k$ and $g_{k}=1$, while $g'_j=1$ 
for all $i\in \{1,\ldots,n\}$. As in the 
previous paragraph it follows that $((g_1,\ldots,g_n),h)$ 
and $((g'_1,\ldots,g'_n),h')$ are adjacent in $\Gamma_{G\wr H}$,
contradicting the fact that they both belong to $S$. This shows 
that $S_H=H_k$, and so $S$ is the stabilizer 
of the point $(w,k)$, completing 
the proof of Proposition~\ref{prop:generalizaion of AhmadiM}.
\end{proof}

\begin{remark}\label{rem:Ahmadi Meagher mistake S_m wr S_n}{\rm
\cite[Proposition 6.5]{AhmadiM}, which claims that  groups $S_m \wr S_n \; (m,n \ge 1)$ have the strict-EKR-property, is wrong, as Corollary~\ref{cor:S_3 wr S_2} shows that $S_3\wr S_2$ does not have the strict-EKR-property. The proof of \cite[Proposition 6.5]{AhmadiM} uses the result that internal direct product $S_m\times S_m \times \cdots \times S_m$ has the strict-EKR-property. However, this is not true for $m \in \{2,3\}$ hence the proof of \cite[Proposition 6.5]{AhmadiM} is correct only for $m\geq 4$.}
\end{remark}

We now investigate when is it that the group $S_m \wr S_n$ has the strict-EKR-property. As explained in Remark \ref{rem:Ahmadi Meagher mistake S_m wr S_n}, one only needs to consider the cases  $m\in \{2,3\}$. We will need the following two lemmas. Observe that $S_n$ has the strict-EKR-property and admits an element fixing exactly one point.
\begin{lemma}\label{lem:S_2 wr H}
Let $H\leq S_n$ be a transitive permutation group having the strict-EKR-property admitting an element fixing exactly one point. Then the group $S_2 \wr H$ has the strict-EKR-property.
\end{lemma}
\begin{proof}
Since $H$ has the strict-EKR-property, it follows from Corollary \ref{cor:wreath1} that $S_2\wr H$ has the EKR-property. Let $I$ be a maximum intersecting set in $S_2\wr H$ containing the identity. As $S_2\wr H$ has the EKR-property, we get that $|I|$ is equal to the order of the stabilizer of $(1,1)$ in $S_2\wr H$, which is further equal to $2^{n-1} |H|/n$ by the Orbit-stabilizer lemma. For $h\in H$ let $I_h=I\cap (S_2^n \times \{h\})$ and observe that $|I_h| \le 2^{n-1}$. Let $I_H=\{h \in H\mid I_h\neq \emptyset \}$ and observe that $I_H$ is an intersecting set in $H$. As $H$ has the strict-EKR-propery, we have that $|I_H| \le |H_1| = |H|/n$ by the Orbit-stabilizer lemma. Therefore, we have that
$$
   \frac{2^{n-1} |H|}{n} =|I| = \sum_{h \in I_H} |I_h| \le |I_H| 2^{n-1} \le \frac{2^{n-1} |H|}{n}.
$$
This shows that $|I_h|=2^{n-1}$ for each $h \in I_H$ and that $|I_H|=|H|/n$, and so $I_H$ is an intersecting set in $H$ of maximum size. Hence $I_H$ is a stabilizer of a point. Without loss of generality, we assume that $I_H$ is the stabilizer of $1$ in $H$.

Let $d \in I_H$ be such that $fix(d)=\{1\}$ (such an element exists by the assumption on $H$). Since the identity is contained in the intersecting set $I$, and $|I_d|=2^{n-1}$, it follows that $I_d=(\{1\}\times S_2\times \cdots \times S_2) \times \{d\}$. 

Suppose now that there exists $\mathbf{g}=((g_1,\ldots,g_n),h) \in I$ with $g_1\neq 1$. Let $\mathbf{g}'=((1,g'_2,\ldots,g'_n),d)$ where 
$g'_i=g_it$ with $t=(1\,2)\in S_2$. 
Observe that $\mathbf{g}$ and $\mathbf{g}'$ both belong to $I$, 
which contradicts the fact that they are adjacent in $\Gamma_{S_2\wr H}$. 
This shows that for each $h\in I_H$ we have that $I_h=(\{1\}\times S_2 \times \cdots \times S_2) \times \{h\}$ and  we conclude that $I=(\{1\} \times S_2 \times \cdots \times S_2) \times I_H$. Hence $I$ is 
the stabilizer of the point $(1,1)$.
\end{proof}

\begin{lemma}\label{lem:S_3 wr H}
Let $H\leq S_n$ be a transitive permutation group having the strict-EKR-property admitting an element fixing exactly one point, and let $n\geq 3$.  Then the group $S_3\wr H$ has the strict-EKR-property.
\end{lemma}
\begin{proof}
Since $H$ has the strict-EKR-property, it follows from Corollary \ref{cor:wreath1} that $S_3\wr H$ has the EKR-property. Let $I$ be a maximum intersecting set in $S_3\wr H$ containing the identity. For $h\in H$ let $I_h=I\cap (S_3^n \times \{h\})$ and let  $I_H=\{h \in H\mid I_h\neq \emptyset \}$. Similarly as in the proof of Lemma \ref{lem:S_2 wr H} we get that $|I_h|=2\cdot 6^{n-1}$ for every $h \in I_H$ and that $I_H$ is a stabilizer of a point. Without loss of generality, we may assume that $I_H$ is the stabilizer of $1$ in $H$.

By the assumption on $H$ there exists $d \in I_H$  such that $fix(d)=\{1\}$. 
It follows that for every $(g_1,\ldots,g_n),d) \in I$ we have that $g_1$ fixes a 
point. Moreover, if $g_1\neq 1$, then every element in $I_d$ has the first 
coordinate equal to $1$ or $g_1$, as otherwise we would obtain two 
adjacent elements of $I_d$. 
Observe that not every element in  $I_d$ has the first coordinate 
equal to $1$, as then we would have 
$|I_d|\leq 6^{n-1}$. Let $f\in S_3\setminus \{1\}$ be a 
non-derangement that appears as the first coordinate of some element from 
$I_d$. Then $I_d=(\{1,f\} \times S_3 \times \cdots \times S_3) \times \{d\}$. 

Suppose that there exists $\mathbf{g}=((g_1,\ldots,g_n),h) \in I$ with $g_1\not \in \{1,f\}$. Further, let $\mathbf{g}'=((g'_1,g'_2,\ldots,g'_n),d)\in I_d$ where 
$g'_1=1$ if $g_1\in \{(1\,2\,3),(1\,3\,2)\}$ and $g'_1=f$ otherwise, while $g'_i=g_i(1\,2\,3)$ for $i\neq 1$. Observe that $\mathbf{g}$ and $\mathbf{g}'$ both belong to $I$, contradicting the fact that they are adjacent in $\Gamma_{S_3\wr H}$. This shows that for each $h\in I_H$ we have  $I_h=(\{1,f\}\times S_3 \times \cdots \times S_3) \times \{h\}$. We conclude that $I=(\{1,f\} \times S_3 \times \cdots \times S_3) \times I_H$. Hence $I$ is a stabilizer of the point $(x,1)$, where $fix(f)=\{x\}$.
\end{proof}

\begin{proposition}\label{pro:S_m wr S_n the strict EKR}
Let $m$ and $n$ be positive integers. The group $S_m \wr S_n$ has the strict-EKR-property if and only if ($m,n)\neq (3,2)$. 
\end{proposition}
\begin{proof}
The proof follows combining Proposition~\ref{prop:generalizaion of AhmadiM}, Corollary~\ref{cor:S_3 wr S_2} and Lemmas~\ref{lem:S_2 wr H} and \ref{lem:S_3 wr H}.

\end{proof}


 \end{document}